\documentclass[twoside,11pt]{article}

\usepackage{graphicx, subfig}
\usepackage[top=1in,bottom=1in,left=1in,right=1in]{geometry}
\usepackage[sort&compress]{natbib} 
\usepackage{amsfonts, amsmath, amssymb, amsthm, constants, bbm}
\usepackage{MnSymbol}
\usepackage{enumitem}
\usepackage[mathscr]{eucal}
\usepackage{booktabs}

\usepackage{color}
\definecolor{darkred}{RGB}{100,0,0}
\definecolor{darkgreen}{RGB}{0,100,0}
\definecolor{darkblue}{RGB}{0,0,150}

\usepackage{hyperref}
\hypersetup{colorlinks=true, linkcolor=darkred, citecolor=darkgreen, urlcolor=darkblue}
\usepackage{url}

\newtheorem{thm}{Theorem}

\theoremstyle{remark}

\newtheorem{rem}{Remark}

\theoremstyle{definition}
\newtheorem{con}{Contribution}

\def\beq{\begin{equation}} 
\def\eeq{\end{equation}}
\def\beqn{\begin{eqnarray*}}
\def\eeqn{\end{eqnarray*}}
\def\Bitem{\begin{itemize}\setlength{\itemsep}{.2in}}
\def\bitem{\begin{itemize}\setlength{\itemsep}{.05in}}
\def\eitem{\end{itemize}}
\def\Benum{\begin{enumerate}\setlength{\itemsep}{.2in}}
\def\benum{\begin{enumerate}\setlength{\itemsep}{.05in}}
\def\eenum{\end{enumerate}}
\def\bmult{\begin{multline*}}
\def\emult{\end{multline*}}
\def\bcenter{\begin{center}}
\def\ecenter{\end{center}}
\def\bframe{\begin{frame}}
\def\eframe{\end{frame}}

\newcommand{\thmref}[1]{Theorem~\ref{thm:#1}}

\newcommand{\secref}[1]{Section~\ref{sec:#1}}
\newcommand{\figref}[1]{Figure~\ref{fig:#1}}

\def\cI{\mathcal{I}}
\def\cJ{\mathcal{J}}

\def\cN{\mathcal{N}}

\def\cS{\mathcal{S}}

\def\bR{\mathbf{R}}

\def\bX{\mathbf{X}}

\def\bx{\mathbf{x}}

\def\bbI{\mathbb{I}}

\def\bbR{\mathbb{R}}
\def\bbS{\mathbb{S}}

\newcommand{\E}{\operatorname{\mathbb{E}}}
\renewcommand{\P}{\operatorname{\mathbb{P}}}
\newcommand{\Var}{\operatorname{Var}}

\def\eps{\varepsilon}

\def\1{\mathbbm{1}}
\newcommand{\IND}[1]{\bbI\{ #1 \}}

\def\A{\mathscr{A}}
\def\B{\mathscr{B}}
\def\X{\bX}
\def\R{\bR}
\def\scan{\textsc{scan}}

\definecolor{purple}{rgb}{0.4,.1,.9}

\urldef\eacurl\url{http://www.math.ucsd.edu/~eariasca/}
\urldef\ylurl\url{http://www.math.ucsd.edu/~yul085/}

\begin{document}
\thispagestyle{empty}

\title{Distribution-free Detection of a Submatrix}
\author{
Ery Arias-Castro\thanks{University of California, San Diego --- \eacurl} 
\and 
Yuchao Liu\thanks{University of California, San Diego --- \ylurl}
}
\date{}
\maketitle

\vspace{-0.3in}
\bcenter\small
\begin{tabular}{p{0.8\textwidth}}
\toprule
{\em Abstract.}
We consider the problem of detecting the presence of a submatrix with larger-than-usual values in a large data matrix.  This problem was considered in \citep{butucea2013detection} under a one-parameter exponential family, and one of the test they analyzed is the scan test.  
Taking a nonparametric stance, we show that a calibration by permutation leads to the same (first-order) asymptotic performance.  This is true for the two types of permutations we consider.  We also study the corresponding rank-based variants and precisely quantify the loss in asymptotic power.
\\  
\bottomrule
\end{tabular}
\ecenter


\section{Introduction} \label{sec:intro}

Biclustering has emerged as an important set of tools in bioinformatics, in particular, in the analysis of gene expression data \citep{cheng2000biclustering}.  It comes in different forms, and in fact, the various methods proposed under that umbrella may target different goals.  See \citep{madeira2004biclustering} for a survey.
Here we follow \citep{shabalin2009finding}, where the problem is posed as that of discovering a submatrix of unusually large values in a (large) data matrix.  For example, in the context of a microarray dataset, the data matrix is organized by genes (rows) and samples (columns).  We let $\X = (X_{ij})$ denote the matrix, $M$ denote the number of rows and $N$ denote the number of columns, so the data matrix $\X$ is $M$-by-$N$.  

\subsection{Submatrix detection} \label{sec:intro-detection}
In its simplest form, there is only one submatrix to be discovered.  In that context, the detection problem is that of merely detecting of the presence of an anomalous (or unusual) submatrix, which leads to a hypothesis testing problem.  This was considered in \citep{butucea2013detection} from a minimax perspective.  Their work relies on parametric assumptions.  For example, in the normal model, they assume that the $X_{ij}$'s are independent and normal, with mean $\theta_{ij}$ and unit variance.  Under the null hypothesis, $\theta_{ij} = 0$ for all $i \in [M] := \{1, \dots, M\}$ and all $j \in [N]$.  Under the alternative, there is a $m$-by-$n$ submatrix indexed by $\cI_{\rm true} \subset [M]$ and $\cJ_{\rm true} \subset [N]$ such that 
\beq\label{theta}
\theta_{ij} \ge \theta_\ddag, \quad \forall (i,j) \in \cI_{\rm true} \times \cJ_{\rm true},
\eeq
while $\theta_{ij} = 0$ otherwise.  Here $\theta_\ddag > 0$ controls the signal-to-noise ratio.  In that paper, Butucea and Ingster precisely establish how large $\theta_\ddag$ needs to be as a function of $(M,N,m,n)$ in order for there to exist a procedure that has (worst-case) risk tending to zero in the large-sample limit (i.e., as the size of the matrix grows).  They consider two tests which together are shown to be minimax optimal.
One is the sum test based on
\beq\label{sum}
\textsc{sum}(\X) = \sum_{i \in [M]} \sum_{j \in [N]} X_{ij}.
\eeq
It is most useful when the submatrix is large.
The other one is the scan test which, when the submatrix size is known (meaning $m$ and $n$ are known) is based on
\beq\label{scan}
\scan(\X)= \max_{\cI \subset [M], |\cI| = m} \quad \max_{\cJ \subset [N], |\cJ| = n} \quad \sum_{i \in \cI} \sum_{j \in \cJ} X_{ij}.
\eeq
When $m$ and $n$ are unknown, one can perform a scan test for each $(m,n)$ in some range of interest and control for multiple testing using the Bonferroni method.  From \citep{butucea2013detection}, and also from our own prior work, we know that the resulting procedure achieves the same first-order asymptotic performance.

To avoid making parametric assumptions, some works such as \citep{barry2005significance,hastie2000gene} have suggested a calibration by permutation.  
We consider two somewhat stylized permutation approaches:
\bitem
\item {\em Unidimensional permutation.}
The entries are permuted within their row.  (One could permute within columns, which is the same after transposition.)  
\item {\em Bidimensional permutation.}
The matrix is vectorized, the entries are permuted uniformly at random as one would in a vector, and the matrix is reformed.  \eitem
The first method is most relevant when one is not willing to assume that the entries in different rows are comparable.  It is appealing in the context of microarray data and was suggested, for example, in \citep{hastie2000gene}.
The second method is most relevant in a setting where all the variables are comparable.  
In the parlance of hypothesis testing, the first method derives from a model where the entries within each row are exchangeable under the null, while the second method arises when assuming that all the entries are exchangeable under the null.

\begin{con}[Calibration by permutation]
We analyze the performance of the scan test when calibrated using one of these two permutation approaches.  We show that, regardless of the variant,  the resulting test is (first order) asymptotically as powerful as a calibration by Monte Carlo with full knowledge of the parametric model.  We prove this under some standard parametric models.  
\end{con}

\begin{rem}
We focus on the scan statistic \eqref{scan} and abandon the sum statistic \eqref{sum} for at least two reasons: 1) the sum statistic cannot be calibrated without knowledge of the null distribution; 2) the sum statistic is able to surpass the scan statistic when it is impossible to locate the submatrix with any reasonable accuracy, which is somewhat less interesting to the practitioner.
\end{rem}

\medskip
A calibration by permutation is computationally intensive in that it requires the repeated computation of the test statistic on permuted data.  In practice, several hundred permutations are used, which can cause the method to be rather time-consuming.  A possible way to avoid this is to use ranks, which was traditionally important before the availability of computers with enough computational power.  \citep{MR758442} is a classical reference.  In line with the two permutation methods described above, we consider the corresponding methods for ranking the entries:
\bitem
\item {\em Unidimensional ranks.}
The entries are ranked relative to the other entries in their row.  
\item {\em Bidimensional ranks.}
The entries are ranked relative to the all other entries. 
\eitem
The use of ranks has the benefit of only requiring calibration (typically done on a computer nowadays) once for each matrix size $M \times N$.  It has the added benefit of yielding a method that is much more robust to outliers.

\begin{con}[Rank-based method]
We analyze the performance of the scan test when the entries are replaced by their ranks following one of the two methods just described.  We show that, regardless of the variant, there is a mild loss of asymptotic power, which we precisely quantify.  We do this under some standard parametric models.  
\end{con}

\subsection{More related work}
The scan statistic \eqref{scan} is computationally intractable and there has been efforts to offer alternative approaches.  We already mentioned \citep{shabalin2009finding}, which proposes an alternate optimization strategy: given a set of rows, optimize over the set of columns, and vice versa, alternating in this fashion until convergence to a local maximum.  This is the algorithm we use in our simulations.  It does not come with theoretical guarantees (other than converging to a local maximum) but performs well numerically.
A spectral method is proposed in \citep{cai2015computational} and a semidefinite relaxation is proposed in \citep{chen2014statistical}.  These methods can be run in time polynomial in the problem size (meaning in $M$ and $N$).  \citep{ma2015computational} establishes a lower bound based on the Planted Clique Problem that strictly separates the performance of methods that run in polynomial time from the performance of the scan statistic.   

Our work here is not on the computational complexity of the problem.  Rather we assume that we can compute the scan statistic and proceed to study it.  In effect, we contribute here to a long line of work that studies permutation and rank-based methods for nonparametric inference.  Most notably, we continue our recent work \citep{arias2015distribution} where we study the detection problem under a similar premise but under much more stringent structural assumptions.  The setting there would correspond to an instance where the submatrix is in fact a block, meaning, that $\cI_{\rm true}$ and $\cJ_{\rm true}$ are of the form $\cI_{\rm true} = \{i+1, \dots, i+k\}$ and $\cJ_{\rm true} = \{j+1, \dots, j+l\}$.  The present setting assumes much less structure.  The related applications are very different in the end.  Nevertheless, the technical arguments developed there apply here with only minor adaptation.
The main differences are that we consider two types of permutation and ranking protocols. 

\subsection{Content}
The rest of the paper is organized as follows.
In \secref{parametric} we describe a parametric setting where likelihood methods have been shown to perform well.  This parametric setting will serve as benchmark for the nonparametric methods that ensue.
In \secref{permutation} we consider the detection problem and study the scan statistic with each of the two types of calibration by permutation.
In \secref{rank-detect} we consider the same problem and study the rank-based scan statistic using each of the two types of rankings.
In \secref{numerics} we present some numerical experiments on simulated data.
All the proofs are in \secref{proofs}.

\section{The parametric scan}
\label{sec:parametric}

Following the classical line in the literature on nonparametric tests, we will evaluate the nonparametric methods introduced later on a family on parametric models.  As in \citep{butucea2013detection}, and in our preceding work \citep{arias2015distribution}, we consider a one-parameter exponential family in natural form.

To define such a family, fix a probability distribution $\nu$ on the real line with zero mean and unit variance, and with a sub-exponential right tail, specifically meaning that $\varphi(\theta) := \int_\bbR e^{\theta x} \nu({\rm d}x) < \infty$ for some $\theta > 0$.  
Let $\theta_\star$ denote the supremum of all such $\theta>0$.  (Note that $\theta_\star$ may be infinite.)
The family is then parameterized by $\theta \in [0, \theta_\star)$ and has density with respect to $\nu$ defined as 
\begin{align} \label{model}
f_\theta(x)  = \exp \{\theta x - \log \varphi(\theta)\}.
\end{align}
By varying $\nu$, we obtain the normal (location) family, the Poisson family (translated to have zero mean), and the Rademacher family.

Such a parametric model is attractive as a benchmark because it includes these popular models and also because likelihood methods are known to be asymptotically optimal under such a model.  \cite{butucea2013detection} showed this to be the case for the problem of detection, where  the generalized likelihood ratio test is based on the scan statistic \eqref{scan}.

Under such a parametric model, the detection problem is formalized as a hypothesis testing problem where $\nu$ plays the role of null distribution.  In detail, suppose that the submatrix is known to be $m \times n$.  The search space is therefore 
\beq
\bbS_{m,n} := \big\{\cS = \cI \times \cJ : \cI \subset [M], |\cI| = m \text{ and } \cJ \subset [N], |\cJ| = n\big\}.
\eeq
We assume that the $X_{ij}$'s are independent with $X_{ij} \sim f_{\theta_{ij}}$, and the testing problem is 
\beq
H_0: \theta_{ij} = 0, \quad \forall (i,j) \in [M] \times [N],
\eeq
versus
\beq
H_1: \exists \cS_{\rm true} \in \bbS_{m,n} \text{ such that } 
\begin{cases} 
\theta_{ij} \ge \theta_\ddag, & \forall (i,j) \in \cS_{\rm true}, \\
\theta_{ij} = 0, & \text{otherwise}.
\end{cases}
\eeq
Here $\theta_\ddag$ controls the signal-to-noise ratio is assumed to be known in this formulation.

In this context, we have the following.

\begin{thm}[\cite{butucea2013detection}] \label{thm:butucea}
Consider an exponential model as described above, with $\nu$ having finite fourth moment.  Assume that 
\beq\label{basic}
M, N, m, n \to \infty, \quad \frac{m}M, \frac{n}N \to 0, \quad \frac{\log (M \vee N)}{m \wedge n} \to 0.
\eeq
Then the sum test based on \eqref{sum}, at any fixed level $\alpha > 0$, has limiting power~1 when 
\beq
\theta_\ddag \frac{m n}{\sqrt{MN}} \to \infty.
\eeq
Then the scan test based on  \eqref{scan}, at any fixed level $\alpha > 0$, has limiting power~1 when 
\beq \label{scan-cond}
\liminf \frac{\theta_\ddag \sqrt{mn}}{\sqrt{2 (m\log \frac{M}{m} + n\log \frac{N}{n})}} > 1.
\eeq
Conversely, the following matching lower bound holds.  Assume in addition that $\log M \asymp \log N$ and $m \asymp n$.  
Then any test at any fixed level $\alpha > 0$ has limiting power at most $\alpha$ when
\beq
\theta_\ddag \frac{m n}{\sqrt{MN}} \to 0 \quad \text{and} \quad \liminf \frac{\theta_\ddag \sqrt{mn}}{\sqrt{2 (m\log \frac{M}{m} + n\log \frac{N}{n})}} < 1.
\eeq
\end{thm}
We note that \cite{butucea2013detection} derived their lower bound under slightly weaker assumptions on $M,N,m,n$.  

\begin{rem}
Proper calibration in this context is based on knowledge of the null distribution $\nu$.  In more detail, consider a test that rejects for large values of a statistic $T(\X)$.  Assuming a desired level of $\alpha > 0$ and that $\nu$ is either diffuse or discrete (for simplicity), the critical value for $T$ is set at $t_\alpha$, where $t_\alpha = \inf\{t :  \nu(T(\X) \ge t) \le \alpha\}$.  The test is then $\IND{T(\X) \ge t_\alpha}$.
In practice, $t_\alpha$ may be approximated by Monte Carlo sampling.
\end{rem}

\section{Permutation scan tests}
\label{sec:permutation}

In the previous section we described the work of \cite{butucea2013detection}, who in certain parametric models show that the sum test \eqref{sum} and scan test \eqref{scan} are jointly optimal for the problem of detecting a submatrix.  This is so if they are both calibrated with full knowledge of the null distribution (denoted $\nu$ earlier).

What if the null distribution is unknown?  
A proven approach is via permutation.  This is shown to be optimal in some classical settings \citep{MR2135927} and was recently shown to also be optimal in more structured detection settings \citep{arias2015distribution}.  We prove that this is also the case in the present setting of detecting a submatrix.
We consider the two types of permutation, unidimensional and bidimensional, described in \secref{intro-detection}.
More elaborate permutation schemes have been suggested, e.g., in \citep{barry2005significance}, but these are not considered here, in part to keep the exposition simple.  Indeed, we simply aim at showing that a calibration by permutation performs very well in the present context.

Let $\Pi$ be a subgroup of permutations of $[M] \times [N]$, identified with $[MN]$.  Then a calibration by permutation of the scan statistic (or any other statistic) yields the P-value
\beq \label{permscan}
\mathfrak{P}(\X) = \frac{\# \{ \pi \in \Pi : \scan(\X_\pi) \geq \scan(\X) \}}{|\Pi|} ,
\eeq
where $\X_\pi = (X_{\pi(i,j)})$ is the matrix permuted by $\pi$.   
The permutation scan test at level $\alpha$ is the test $\IND{\mathfrak{P}(\X) \le \alpha}$.
It is well-known that this this a valid P-value, in the sense that, under the null, it dominates the uniform distribution on $[0,1]$ \citep{MR2135927}.  (This remains true of a Monte Carlo approximation.)

The set of unidimensional permutations, denoted $\Pi_1$, is that of all permutations that permute within each row, while the set of bidimensional permutations, denoted $\Pi_2$, is simply the set of all permutations.
Obviously, $\Pi_1 \subset \Pi_2$ with $|\Pi_1| = (N!)^M$ and $|\Pi_2| = (MN)!$, and they are both groups.\footnote{The group structure is important.  See the detailed discussion in \citep{hemerik2014exact}.}

\begin{thm} \label{thm:permutation}
Consider an exponential model as described in \secref{parametric}.  In addition to \eqref{basic}, assume
\beq\label{basic1}
\log^3(M \vee N)/(m \wedge n) \to 0,
\eeq 
and that either {\em (i)} $\nu$ has support bounded from above, or {\em (ii)} $\max_{i,j} \theta_{ij} \le \bar\theta$ for some $\bar\theta < \theta_\star$ fixed. 
Let the group of permutations $\Pi$ be either $\Pi_1$ or $\Pi_2$; if $\Pi = \Pi_1$, we require that $\varphi(\theta) < \infty$ for some $\theta < 0$. Then the permutation scan test based on \eqref{permscan}, at any fixed level $\alpha > 0$, has limiting power~1 when \eqref{scan-cond} holds.
\end{thm}

The additional condition (on $\nu$ or the nonzero $\theta_{ij}$'s) seems artificial, but just as in \citep{arias2015distribution}, we are not able to eliminate it.  Other than that, in view of \thmref{butucea} we see that the permutation scan test --- just like the parametric scan test --- is optimal to first-order under a general one-parameter exponential model.

\section{Rank-based scan tests}
\label{sec:rank-detect}

Rank tests are classical special cases of permutation tests \citep{MR758442}.  Traditionally, when computers were not as readily available and not as powerful, permutation tests were not practical, but rank tests could still be, as long as calibration had been done once for the same (or a comparable) problem size.  Another well-known advantage of rank tests is their robustness to outliers.  

We consider the two ranking protocols described in \secref{intro-detection}.  After the observations are ranked, the distribution under the null is the permutation distribution, either uni- or bi-dimensional depending on the ranking protocol.  This is strictly true under an appropriate exchangeability condition, which holds in the null model we consider here where all observations are IID.  
In fact, the unidimensional rank scan test is a form of unidimensional permutation test, and the bidimensional rank scan is a form of bidimensional permutation test, each time, the statistic being the rank scan
\beq\label{rank-scan}
\scan(\R)= \max_{\cI \subset [M], |\cI| = m} \quad \max_{\cJ \subset [N], |\cJ| = n} \quad \sum_{i \in \cI} \sum_{j \in \cJ} R_{ij},
\eeq
where $\R = (R_{ij})$ is the matrix of ranks.
 
Rank tests have been studied in minute detail in the classical setting \citep{MR758442,hajek1967theory}.  Typically, this is done, again, by comparing their performance with the likelihood ratio test in the context of some parametric model.  Typically, there is some loss in efficiency, unless one tailors the procedure to a particular parametric family.\footnote{Actually, \cite{hajek1962asymptotically} proposes a more complex method that avoids the need for knowing the null distribution.}  
Such a performance analysis was recently carried out for the rank scan in more structured settings \citep{arias2015distribution}.  We again extend this work here and obtain the following. 

Define 
$$
\Upsilon = \E (Z \1_{(Z > Y)})  + \frac{1}{2} \E (Z \1_{(Z = Y)})
$$
where $Y,Z$ are IID with distribution $\nu$.  (This is the same constant introduced by \cite{arias2015distribution}.)

\begin{thm} \label{thm:ranktest}
Consider an exponential model as described in \secref{parametric}.  Assume that \eqref{basic} holds. 
Let the group of permutations $\Pi$ be either $\Pi_1$ or $\Pi_2$.
The rank scan test at any fixed level $\alpha > 0$ has limiting power~1 when
\beq \label{rank-cond}
\liminf \frac{\theta_\ddag \sqrt{mn}}{\sqrt{2 (m\log \frac{M}{m} + n\log \frac{N}{n})}} > \frac{1}{2 \sqrt{3} \Upsilon}.
\eeq
\end{thm}

The proof is omitted as it is entirely based on an adaptation of that of \thmref{permutation} and arguments given in \citep{arias2015distribution} to handle the rank moments.  

Compared to the (optimal) performance of the parametric and permutation scan tests in the same setting (\thmref{butucea} and \thmref{permutation}), we see that there is a loss in power.  However, the loss can be quite small.  For example, as argued in \citep{arias2015distribution}, in the normal model $1/(2 \sqrt{3} \Upsilon) \le \sqrt{\pi/3} \approx 1.023$.

\section{Numerical experiments}
\label{sec:numerics}

We performed some numerical experiments\footnote{In the spirit of reproducible research, our code is publicly available at \url{https://github.com/nozoeli/NPDetect}} to assess the accuracy of our asymptotic theory.
To do so, we had to deal with two major issues in terms of computational complexity. 
The first issue is the computation of the scan statistic defined in \eqref{scan}.  There are no known computationally tractable method for doing so.  As \cite{butucea2013detection} did, we opted instead for an approximation in the form of the alternate optimization (or hill-climbing) algorithm of \cite{shabalin2009finding}.  Since in principle this algorithm only converges to a local maximum, we run the algorithm on several random initializations and take the largest output. 
The second issue is that of computing the permutation P-value defined in \eqref{permscan}.  (This is true for the permutation test and also for the special case of the rank test.)  Indeed, examining all possible permutations in $\Pi$ (either $\Pi_1$ or $\Pi_2$) is only feasible for very small matrices.  As usual, we opted for Monte Carlo sampling.  Specifically, we picked $\pi_1, \ldots \pi_B$ IID uniform from $\Pi$ with $B = 500$ in our setup.  We then estimate the permutation P-value by
\beq
\hat{\mathfrak{P}}(\X) = \frac{\# \{b \in [B] : \scan(\X_{\pi_b}) \geq \scan(\X) \}}{B+1}.
\eeq

We mention that when rank methods are applied, the ties in the data are broken randomly.  

\paragraph{Simulation setup}
Our simulation strategy is as follows. 
A data matrix $\X$ of size $M \times N$ is generated with the anomaly as $[m] \times [n]$.  All the entries of $\X$ are independent with distribution $f_0$ (same as $\nu$) except for the anomalous ones which have distribution $f_{\theta_\ddag}$ for some $\theta_\ddag > 0$.  
We compare the permutation tests and rank tests (unidimensional and bidimensional) with the scan test calibrated by Monte Carlo (using 500 samples), which serves as an oracle benchmark as it has full knowledge of the null distribution $f_0$.  By construction, all tests have the prescribed level.  
As we increase $\theta_\ddag$, the P-values of the different tests are recorded.
Each setting is repeated 200 times.

As one of the main purposes of our simulations is to confirm our theory, we zoom in on the region near the critical value 
\beq
\theta_{\text{crit}} = \sqrt{\frac{2 (m \log \frac{M}{m} + n \log \frac{N}{n})}{mn}},
\eeq
which comes from \eqref{scan-cond}.
Specifically, we increase $\theta_\ddag$ from $0.5 \times \theta_{\text{crit}} $ to $1.5\times \theta_{\text{crit}}$ with step size $0.125\times  \theta_{\text{crit}}$ to explore the behavior of P-values around the critical value.

\begin{figure}%
	\centering
	\subfloat{{\includegraphics[width=\textwidth]{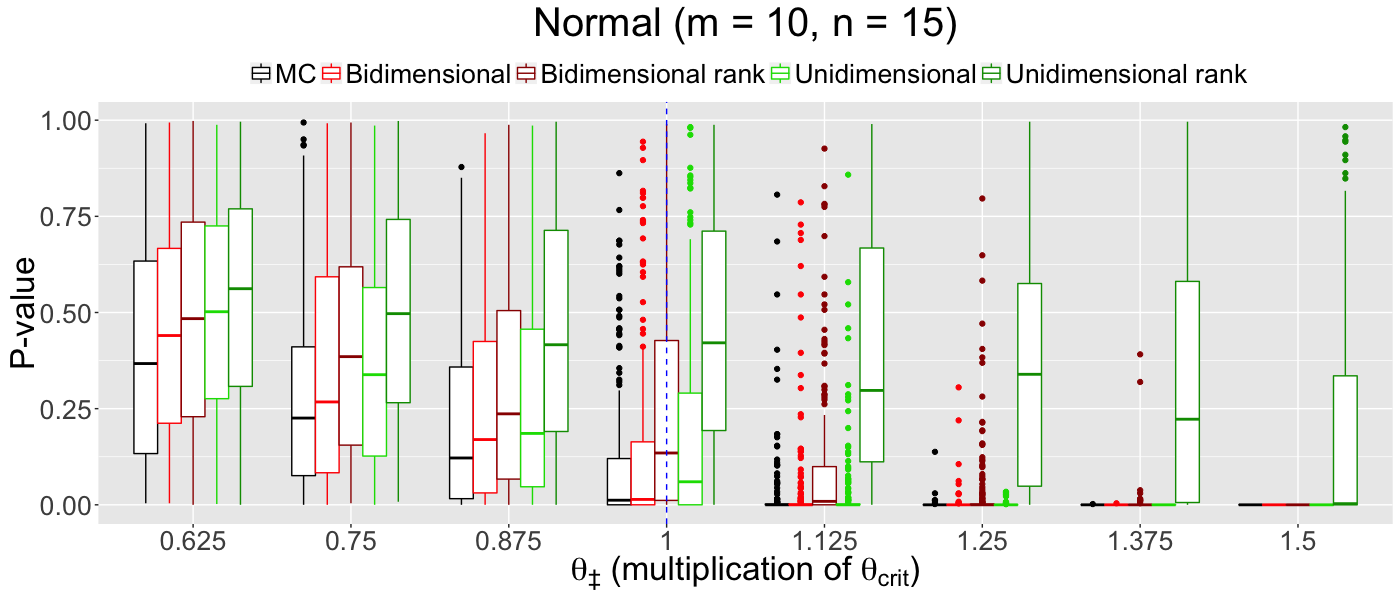} }}\\%
	\subfloat{{\includegraphics[width=\textwidth]{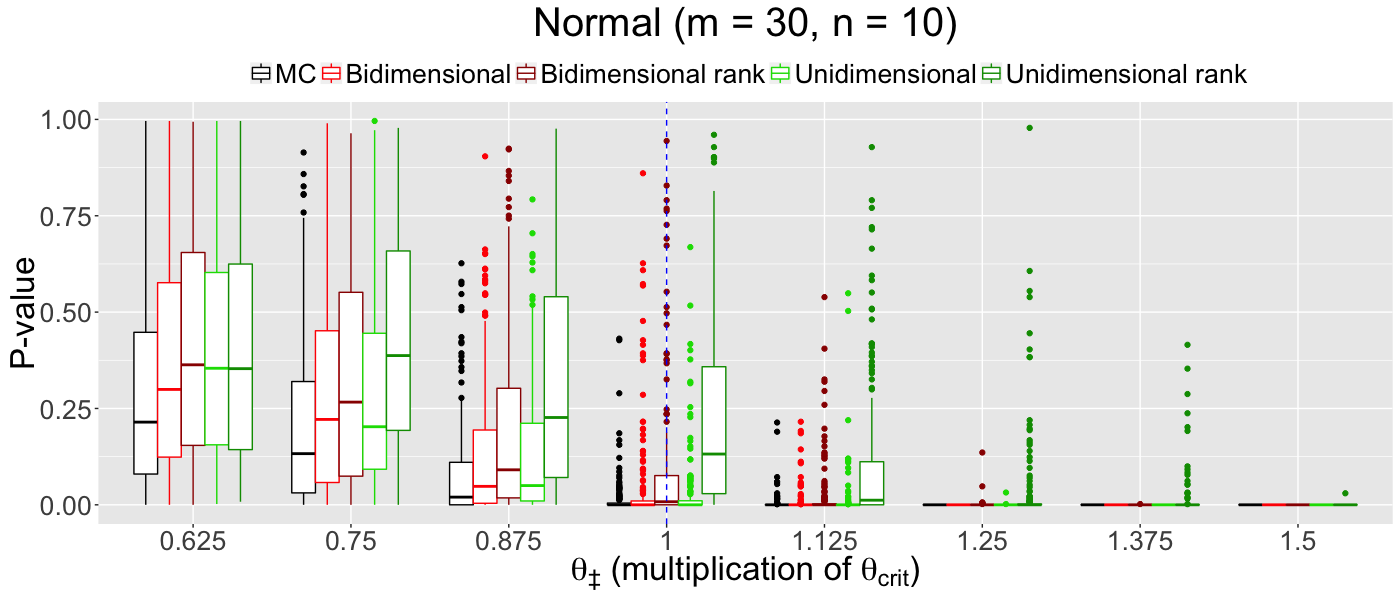} }}%
	\caption{P-values of various forms of scan tests in the normal model}%
	\label{fig:norm}%
\end{figure}

\paragraph{The Normal Case}
Here we generate data from normal family, where $f_\theta$ is $\cN(\theta, 1)$.  
We used two setups, $(M,N,m,n)= (200,100,10,15)$ and $(M,N,m,n)= (200,100,30,10)$, to assess the performance of the tests under different anomaly sizes.  
The resulting boxplots of the averaged P-values are shown in \figref{norm}.

From the plots we see that the P-values are generally very close to $0$ when $\theta_\ddag$ exceeds $\theta_{\text{crit}}$.  When $(m,n) = (10,15)$ the convergence towards 0 is slower, which may be due to the small size of the anomalous submatrix.  As expected, the (oracle) Monte Carlo test is best, followed by the bidimensional permutation test, followed by the unidimensional permutation test.  That said, the differences appear to be minor, which confirms our theoretical findings.

For the rank tests, we observe a similar behavior of the P-values, with the bidimensional showing superiority over the unidimensional rank test, but the loss of power with respect to the oracle test is a bit more substantial, as predicted by the theory.  As shown before, $1/(2 \sqrt 3 \Upsilon) \approx 1.03$ for the standard normal, so that we should place the critical threshold approximately at $1.03 \times \theta_{\text{crit}}$.  This appears to be confirmed in the setting where $(m,n) = (30,10)$.  While the P-values for the rank tests converge relatively slowly when $(m,n) = (10,15)$ (for unidimensional rank test the P-value is close to 0 at $\theta = 1.5 \times\theta_{\text{crit}}$), this may be due to the relatively small size of the anomaly.

\begin{figure}%
	\centering
	\subfloat{{\includegraphics[width=\textwidth]{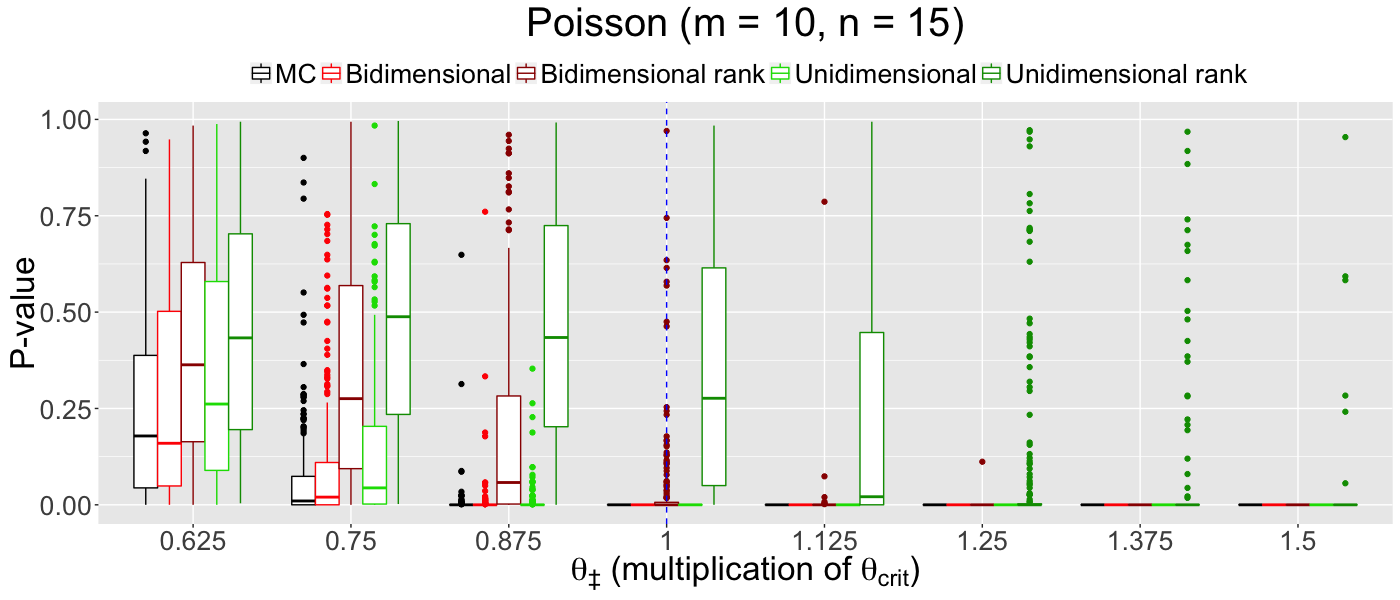} }}\\%
	\subfloat{{\includegraphics[width=\textwidth]{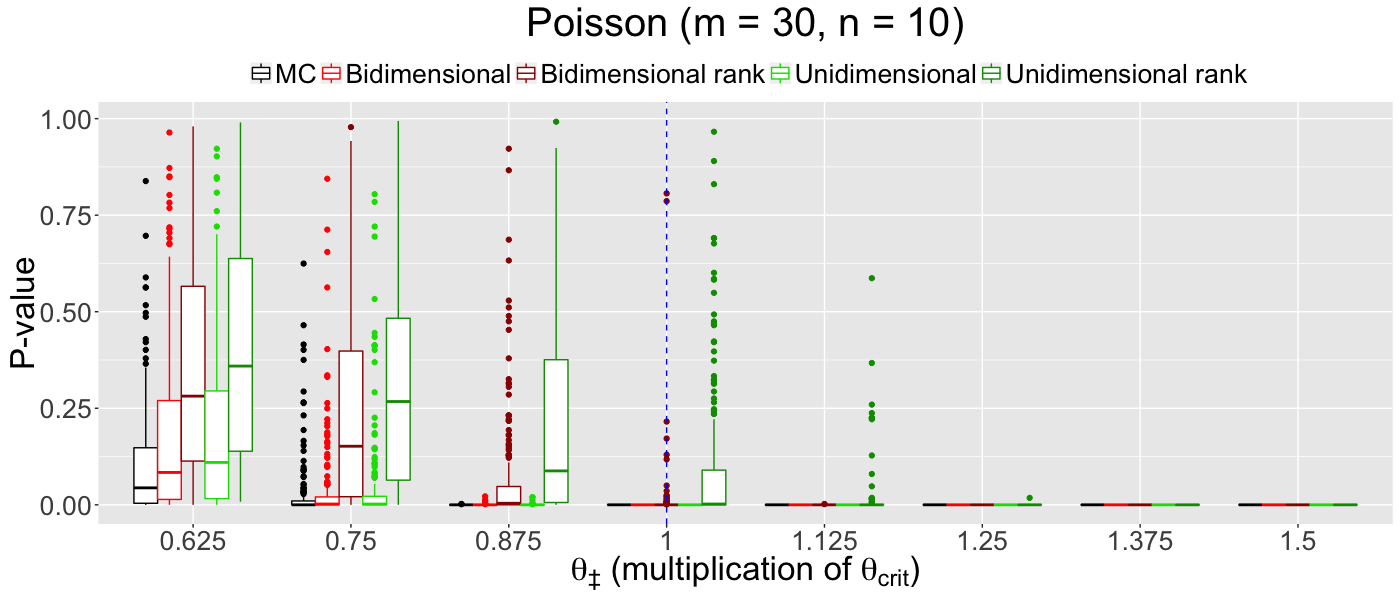} }}%
	\caption{P-values of various forms of scan tests in the Poisson model}%
	\label{fig:poi}%
\end{figure}

\paragraph{The Poisson Case}
As another example, we consider the Poisson family, where $f_\theta$ corresponds to ${\rm Poisson}(e^\theta) - 1$. The data matrix and anomaly sizes are the same as they are in the normal case. The resulting boxplots of the P-values are shown in \figref{poi}.  Overall, we observe a similar behavior of the P-values.

\section{Proofs}
\label{sec:proofs}

\subsection{Preliminaries}

We start with some preliminary results that already appear, one or another, in our previous work \citep{arias2015distribution}.
First, for any one-parameter exponential family $(f_\theta : \theta \in \Theta)$ with a standardized base distribution $\nu$, as we consider to be here,
\beq\label{mean-exp}
\E_\theta(X) \ge \theta, \quad \forall \theta \in \Theta.
\eeq
Next, in the same context, if $\sup\Theta > 0$ (which we assume throughout), then $f_\theta$ has a sub-exponential right tail, which is uniform in $\theta \in \bar\theta$ if $\bar\theta \in \Theta$.
In particular, there is $\bar\gamma$ that depends on $\bar\theta > 0$ such that, if $X_1, \dots, X_k$ are independent, with $X_j \sim f_{\theta_j}$ with $\theta_j \le \bar\theta$, then
\beq\label{upper}
\max_{j \in [k]} X_j \le \bar\gamma \log k, \quad \text{with probability tending 1 as } k \to \infty.
\eeq
By symmetry, if $\inf\Theta < 0$ (which we assume in the case of unidimensional permutations), the same is true on the left.  In particular, $\nu$ itself (corresponding to $\theta = 0$) has a sub-exponential left tail in this case, and this is all that will be used below.  
In particular, there is a constant $\gamma_0 > 0$ such that, if $X_1, \dots, X_k$ are IID $\nu$, then 
\beq\label{lower}
\min_{j \in [k]} X_j \ge - \gamma_0 \log k, \quad \text{with probability tending 1 as } k \to \infty.
\eeq

\subsection{Proof of \thmref{permutation}}
\label{sec:proof-thm-permutation}

The proof arguments are parallel to those of \cite{arias2015distribution}, derived in the context of more structured settings.  For that reason, we only detail the proof of \thmref{permutation} in the case of unidimensional permutations, which, compared to bidimensional permutations, is a bit more different from the setting considered in \citep{arias2015distribution} and requires additional arguments.  Therefore, in what follows, we take $\Pi = \Pi_1$.  Recall that in this case we assume in addition that $\varphi(\theta) < \infty$ for some $\theta < 0$.  
This implies that $\nu$ as sub-exponential tails

\paragraph{Case {\em (i)}}
We first focus on the condition where $\nu$ has support bounded from above and let $b_0$ denote such an upper bound.  (Necessarily, $b_0 > 0$.)  Thus, regardless of the $\theta_{ij}$'s, 
\beq\label{b0}
\P(\textstyle\max_{i,j} X_{ij} \le b_0) = 1.
\eeq
The permutation scan test has limiting power~1 if and only if $\P(\mathfrak{P}(\bX) \le \alpha) \to 1$ under the alternative.  
We show that by proving the stronger claim that $\mathfrak{P}(\bX) \to 0$ in probability under the alternative.

We first work conditional on $\X = \bx$, where $\bx = (x_{ij})$ denotes a realization of $\X = (X_{ij})$.  We may equivalently center the rows of $\X$ before scanning, and the resulting test remains unchanged.  Therefore, we may assume that all the rows of $\bx$ sum to~0.
Let $\zeta = \scan(\bx)$ for short.
We have
\beq
\mathfrak{P}(\bx) = \P(\scan(\bx_\pi) \ge \zeta),
\eeq
where the randomness comes solely from $\pi$, uniformly drawn from $\Pi$.  Using the union bound, we get
\beq
\mathfrak{P}(\bx) \le |\bbS_{m,n}| \max_{\cS \in \bbS_{m,n}} \P\big(\textstyle\sum_{(i,j) \in \cS} x_{\pi(i,j)} \ge \zeta\big).
\eeq
For each $i \in [N]$, let $(A_{ij} : j \in [n])$ be a sample from $(x_{ij} : j \in [N])$ {\em without} replacement and let $A_i = \sum_{j \in [n]} A_{ij}$.  Note that $A_1, \dots, A_M$ are independent and, for $\cS = \cI \times \cJ$, we have
\beq
\sum_{(i,j) \in \cS} x_{\pi(i,j)} \sim \sum_{i \in \cI} A_i.
\eeq   
Fix $\cI \subset [M]$ of size $m$.
Using Markov's inequality and the independence of the $A_i$'s, we get
\beq
\P\big(\textstyle\sum_{i \in \cI} A_i \ge \zeta\big) 
\le e^{-c \zeta} \prod_{i \in \cI} \phi_i(c),
\eeq
where $\phi_i$ is the moment generating function of $A_i$.  
The key is \cite[Th~4]{hoeffding}, which implies that $\phi_i \le \psi_i$, where $\psi_i$ is the moment generating function of $B_i$, where $B_i = \sum_{j \in [n]} B_{ij}$ and $(B_{ij} : j \in [n])$ is a sample from $(x_{ij} : j \in [N])$ {\em with} replacement, meaning that these are IID random variables uniformly distributed in $(x_{ij} : j \in [N])$.
By \eqref{b0}, we have $B_{ij} \le b_0$ with probability one, and the usual arguments leading to the (one-sided) Bernstein's inequality yield the usual bound
\beq
\psi_i(c) \le \exp\left[\frac{n c^2 \sigma_i^2}2  \frac{e^{c b_0} - 1 - c b_0}{c^2 b_0^2/2}\right],
\eeq
where $\sigma_i^2$ is the variance of $B_{i1}$, meaning, $\sigma_i^2 = \frac1N \sum_{j \in [N]} (x_{ij} - \bar x_i)^2$, with $\bar x_i = \frac1N \sum_{j \in [N]} x_{ij}$ being the mean.  Letting $\sigma^2 = \max_{i \in [M]} \sigma_i^2$, we derive 
\beq
e^{-c \zeta} \prod_{i \in \cI} \phi_i(c) 
\le e^{-c \zeta} \prod_{i \in \cI} \exp\left[\frac{n c^2 \sigma_i^2}2  \frac{e^{c b_0} - 1 - c b_0}{c^2 b_0^2/2}\right]
\le e^{-c \zeta} \exp\left[\frac{mn c^2 \sigma^2}2  \frac{e^{c b_0} - 1 - c b_0}{c^2 b_0^2/2}\right],
\eeq
the latter being the usual bound that leads to Bernstein's inequality.  
The same optimization over $c$ yields
\begin{align}
\P\big({\textstyle\sum_{i \in \cI}} A_i \ge \zeta\big) 
&\le \exp\left[- \frac{\zeta^2}{2 mn \sigma^2 + \frac23 b_0 \zeta} \right].
\end{align}
We now emphasize the dependency of $\zeta$ and $\sigma^2$ on $\bx$ by adding $\bx$ as a subscript.  
Noting that this bound is independent of $\cI$ (of size $m$), we get
\beq\label{bound-x}
\mathfrak{P}(\bx) \le |\bbS_{m,n}| \exp\left[- \frac{\zeta_\bx^2}{2 mn \sigma^2_\bx + \frac23 b_0 \zeta_\bx} \right].
\eeq

We now free $\X$ and bound $\zeta_\X$ from below, and $\sigma^2_\bX$ from above.  When doing so, we need to take into account that we assumed the rows summed to~0.  When this is no longer the case, $\zeta_\X$ denotes the scan of $\bX$ after centering all the rows.  Let $\bar X_i$ denote the mean of row $i$.  By definition of the scan in \eqref{scan}, 
\beq
\zeta_\bX \ge \zeta_{\rm true} 
:= \sum_{i \in \cI_{\rm true}} \sum_{j \in \cJ_{\rm true}} (X_{ij} - \bar X_i) 
= (1 - \tfrac nN) \sum_{i \in \cI_{\rm true}} \sum_{j \in \cJ_{\rm true}} X_{ij} - \tfrac{n}N \sum_{i \in \cI_{\rm true}} \sum_{j \notin \cJ_{\rm true}} X_{ij}.
\eeq
For the expectation, by \eqref{basic} and \eqref{mean-exp}, we have 
\beq
\E(\zeta_{\rm true}) \ge (1 - \tfrac nN) \sum_{i \in \cI_{\rm true}} \sum_{j \in \cJ_{\rm true}} \theta_{ij} \ge (1 - o(1)) mn \theta_\ddag.
\eeq
For the variance, we have $\Var(X_{ij}) = 1$ when $(i,j) \notin \cS_{\rm true}$ (since $\nu$ has variance~1) and $\Var(X_{ij}) \le \E(X_{ij}^2) \le b_0^2$ always.  Using this, we derive
\beq
\Var(\zeta_{\rm true}) \le mn b_0^2 + (\tfrac{n}N)^2 m N = mn (b_0^2 + \tfrac{n}N) = O(mn).
\eeq
Because of \eqref{basic} and \eqref{scan-cond}, $\E(\zeta_{\rm true}) \gg \sqrt{\Var(\zeta_{\rm true})}$, and thus by Chebyshev's inequality,
\beq \label{avg-bound}
\zeta_{\rm true} = (1 + o_P(1)) \E(\zeta_{\rm true}) \ge (1 + o_P(1)) mn \theta_\ddag.
\eeq
We now bound $\sigma^2_\bx$. For $i\in \cI_{\rm true}$, we have
\beq
\sigma_i^2(\bX) \le \frac 1N \sum_{j \in [N]} X_{ij}^2  = \frac 1N\sum_{j \in \cJ_{\rm true}} X_{ij}^2 +\frac 1N \sum_{j \notin \cJ_{\rm true}} X_{ij}^2 \le \frac{nb^2_0}{N} + \frac 1N \sum_{j \notin \cJ_{\rm true}} X_{ij}^2.
\eeq
For $i \notin \cI_{\rm true}$,
\beq
\sigma_i^2(\bX) \le \frac 1N \sum_{j \in [N]} X_{ij}^2.
\eeq
Therefore 
\beq\label{sigma-bound}
\sigma^2_\bX \stackrel{\rm sto}{\le} 1 + o(1) + \max_{i \in [M]} \frac{1}{N} \sum_{j \in [N]} T_{ij},
\eeq
where $(T_{ij} : (i,j) \in [M] \times [N])$ are IID with distribution that of $X^2 -1 $ when $X \sim \nu$.
Note that $\E(T_{ij}) = 0$ since $\nu$ has variance 1 and 
\beq
\max_{i,j} T_{ij} \le \bar t := b_0^2 \vee (\gamma_0 \log (MN))^2,
\eeq
by \eqref{b0} and when the following event holds
\beq
\A := \big\{\min_{i,j} X_{ij} \ge -\gamma_0 \log (MN)\big\},
\eeq
which by \eqref{lower} happens with probability tending to 1.  
Let $\P_\A$ be the probability conditional on $\A$ and $\E_\A$ the corresponding expectation.  
Let $\mu_A = \E_\A(T_{ij})$ and $\tau^2_\A = \Var_\A(T_{ij}) < \infty$, because $\nu$ has finite fourth moment.
By Bernstein's inequality, for any $c > \mu_\A$,
\begin{align}
\P_\A\bigg(\frac{1}{N} \sum_{j\in [N]} T_{ij} > c  \bigg) 
&\le \exp \bigg[-\frac{N (c - \mu_\A)^2}{2 \tau_\A^2 + \tfrac 23 \bar t c} \bigg].
\end{align}
Then using a union bound
\beq \label{sigma-max}
\P_\A \bigg(\max_{i \in [M]} \frac{1}{N} \sum_{j\in [N]} T_{ij} > c  \bigg) 
\le M \exp \bigg[-\frac{N (c - \mu_\A)^2}{2 \tau_\A^2 + \tfrac 23 \bar t c} \bigg].
\eeq
Taking logs, noting that $\mu_\A \to 0$ and $\tau^2_\A \to \tau^2 := \Var(T_{ij})$, as well as $\bar t = O(\log (MN))$, and using \eqref{basic} and \eqref{basic1}, we see that the RHS tends to 0 for any $c > 0$ fixed.
Therefore $\max_{i \in [M]}\frac{1}{N} \sum_{j\in [N]} T_{ij} = o_P(1)$ conditional on $\A$, and since $\P(\A) \to 1$, also unconditionally. 
Coming back to \eqref{sigma-bound}, we conclude that 
\beq\label{var-bound}
\sigma^2_\bX = 1 + o_P(1).
\eeq

The upper bound on $\zeta_\bX$ and the lower bound on $\sigma^2_\bX$, combined, imply by monotonicity that
\beq
\frac{\zeta_\X^2}{2 mn \sigma_\X^2 + \frac23 b_0 \zeta_\X} 
\ge (1 + o_P(1)) \frac{mn \theta_\ddag^2}{2 + \frac23 b_0 \theta_\ddag}.
\eeq

We also have $|\bbS_{m,n}| = \binom{M}m \binom{N}n$, so that
\beq\textstyle
\log |\bbS_{m,n}| 
= \log \binom{M}m + \log \binom{N}n 
\le (1+o(1)) \Lambda,
\eeq
with 
\beq
\Lambda := \textstyle\big[m \log \frac{M}m + n \log \frac{N}n\big], 
\eeq
where in the last inequality we used \eqref{basic} and the fact that $\log \binom{K}k \le k \log(K/k) + k$ for all integers $1 \le k \le K$.

Coming back to \eqref{bound-x} and collecting all the bounds in between, we find that
\beq
\log \mathfrak{P}(\bX) \le (1+o(1)) \Lambda - (1 + o_P(1)) \frac{mn \theta_\ddag^2}{2 + \frac23 b_0 \theta_\ddag}.
\eeq
Under \eqref{scan-cond}, there is $\eps>0$ such that, eventually,
\beq\label{theta-ddag-bound}
\theta_\ddag \ge (1+\eps) \sqrt{2\Lambda/(mn)}. 
\eeq
When that's the case, we get
\beq\label{final-bound}
\log \mathfrak{P}(\bX) \le (1+o(1)) \Lambda - (1 + o_P(1)) \frac{(1 + 2\eps) \Lambda}{1 + \frac13 b_0 (1+\eps) \sqrt{2\Lambda/(mn)}}.
\eeq
Noting that $\Lambda/(mn) = o(1)$ and $\Lambda \to \infty$ under \eqref{basic}, we get
\beq
\log \mathfrak{P}(\bX) \le - (1 + o_P(1)) 2\eps \Lambda \to -\infty,
\eeq 
which is what we needed to prove. 

\paragraph{Case {\em (ii)}}
We now consider the case where $\theta_{ij} \leq \bar{\theta}$ for all $(i,j) \in [M] \times [N]$ for some $\bar\theta < \theta_*$. 
Although \eqref{b0} may not hold for any $b_0$, we redefine $b_0 = \bar\gamma \log(MN)$, where $\bar\gamma$ depends on $\bar\theta$, and condition on the event 
\beq
\B := \big\{\max_{i,j} X_{ij} \le b_0\big\},
\eeq
which holds with probability tending to 1 by \eqref{upper}.
The bound \eqref{bound-x} holds unchanged (assuming that $\max_{i,j} x_{ij} \le b_0$).  What is different is how $\zeta_\bX$ and $\sigma^2_\bX$ are handled, now that we conditioned on $\B$.  Let $\P_\B$ and $\E_\B$ denote the probability and expectation conditional on $\B$.

We have
\begin{align}
\E_\B(\zeta_{\rm true}) 
&\ge (1 - \tfrac nN) \sum_{i \in \cI_{\rm true}} \sum_{j \in \cJ_{\rm true}} \E_\B(X_{ij}) - \tfrac{n}N \sum_{i \in \cI_{\rm true}} \sum_{j \notin \cJ_{\rm true}} \E_\B(X_{ij}) \\
&\ge (1 - \tfrac nN) \sum_{i \in \cI_{\rm true}} \sum_{j \in \cJ_{\rm true}} \E(X_{ij} | X_{ij} \le b_0) - \tfrac{n}N \sum_{i \in \cI_{\rm true}} \sum_{j \notin \cJ_{\rm true}} \E(X_{ij} | X_{ij} \le b_0) \\
&\ge (1 + o(1)) mn \theta_\ddag.
\end{align}
In the last inequality, for $j \notin \cJ_{\rm true}$ we used the fact that $\E(X_{ij}) = 0$, which implies that $\E_\B(X_{ij}) \le 0$ in that case.  And for $j \in \cJ_{\rm true}$ we used the fact that $\E(X_{ij} | X_{ij} \le b_0) \to \theta_{ij} \ge \theta_\ddag$ combined with a C\`esaro-type argument. 
On the other hand, in a similar way, we also have 
\beq
\Var_\B(\zeta_{\rm true}) = O(mn b_0^2) = O(mn \log^2(MN)).
\eeq
So we still have $\E_\B(\zeta_{\rm true}) \gg \sqrt{\Var_\B(\zeta_{\rm true})}$, by \eqref{basic} and \eqref{scan-cond}, and in addition \eqref{basic1}.  In particular, \eqref{avg-bound} holds under $\B$.
In very much the same way, one can verify that the same is true of \eqref{var-bound}.

From there we get to \eqref{final-bound} in exactly the same way, conditional on $\B$, and then unconditionally since $\P(\B) \to 1$.  
Then, to conclude, we only need to check that $b_0 \sqrt{\Lambda/(mn)} = o(1)$, which is the case by \eqref{basic1}.

\subsection*{Acknowledgements}
This work was partially supported by a grant from the US Office of Naval Research (N00014-13-1-0257) and a grant from the US National Science Foundation (DMS 1223137).

\bibliographystyle{chicago}
\bibliography{ref}

\end{document}